\documentclass[twoside]{amsart}

\usepackage{amsmath,amssymb,amsbsy,amstext,amsxtra,latexsym}
\usepackage{graphicx}
\usepackage{subfig}
\usepackage{enumerate}
\usepackage[all]{xy}

\newtheorem{theorem}{Theorem}[section]

\newtheorem{proposition}{Proposition}[section]

\theoremstyle{remark}

\newtheorem*{remark}{Remark}

\newcommand{\nc}[1]{\ensuremath \overset{#1}{\circ}}
\newcommand{\ntc}[2]{\ensuremath \underset{#2}{\overset{#1}{\circ}}}

\def\mP{{\mathbb P}}
\def\mQ{{\mathbb Q}}
\def\mC{{\mathbb C}}

\def\mZ{{\mathbb Z}}

\def\b{\bigskip}
\def\m{\medskip}

\begin{document}

\title[A simply connected surface of general type with $p_g=0$ and
       $K^2=4$]{A simply connected surface of general type \\
       with $p_g=0$ and $K^2=4$}

\author{Heesang Park, Jongil Park and Dongsoo Shin}

\address{Department of Mathematical Sciences, Seoul National University,
         San 56-1, Sillim-dong, Gwanak-gu, Seoul 151-747, Korea}

\email{hspark@math.snu.ac.kr}

\address{Department of Mathematical Sciences, Seoul National University,
         San 56-1, Sillim-dong, Gwanak-gu, Seoul 151-747, Korea}

\email{jipark@math.snu.ac.kr}

\address{Department of Mathematical Sciences, Seoul National University,
         San 56-1, Sillim-dong, Gwanak-gu, Seoul 151-747, Korea}

\email{dsshin@math.snu.ac.kr}

\date{September 21, 2008}

\subjclass[2000]{Primary 14J29; Secondary 14J10, 14J17, 53D05}

\keywords{$\mQ$-Gorenstein smoothing, rational blow-down,
          surface of general type}

\begin{abstract}
 As the sequel to~\cite{LP07,PPS},
 we construct a simply connected minimal complex surface of general type
 with $p_g=0$ and $K^2=4$ by using a rational blow-down surgery and
 $\mQ$-Gorenstein smoothing theory.
\end{abstract}

\maketitle

\section{Introduction}

 A rational surface satisfies $p_g=q=0$ and it has Kodaira dimension
 $\kappa=-\infty$. Around 1894 Castelnuovo conjectured that a surface
 with $p_g=q=0$ is rational. However the conjecture was soon shown to
 be false by the examples of Enriques: Non-rational surfaces with
 $p_g=q=0$. Castelnuovo also found other examples of non-rational
 surfaces with $p_g=q=0$. Enriques' examples have Kodaira dimension
 $0$ while Castelnuovo's examples have Kodaira dimension $1$.
 Hence smooth surfaces of general type (i.e. Kodaira dimension $2$)
 with $p_g=q=0$ are very interesting from the point of view of
 the history of surfaces with $p_g=q=0$.

 Nowadays a large number of examples of surfaces of general type with
 $p_g=q=0$ are known due to Godeaux, Campedelli, and so on~\cite{BHPV}.
 However it was only in 1983 that the first example of a
 \emph{simply connected} surface of general type with $p_g=0$ appeared
 -- so-called, Barlow surface~\cite{B}. Barlow surface has $K^2=1$.
 The second examples were discovered just recently.
 Motivated by a result of the second author~\cite{P05},
 Y. Lee and the second author constructed a family of simply connected
 minimal complex surfaces of general type with $p_g=0$ and $K^2=1,2$
 by using rational blow-down surgery and $\mQ$-Gorenstein smoothing
 theory~\cite{LP07}. After this construction,
 the authors constructed a family of simply connected
 minimal complex surfaces of general type with $p_g=0$ and $K^2=3$
 by the similar methods~\cite{PPS}.

 In this paper we extend the results~\cite{LP07,PPS} to the case of $K^2=4$.
 That is, we construct a new simply connected minimal surface of
 general type with $p_g=0$ and $K^2=4$ by using a rational blow-down
 surgery and $\mQ$-Gorenstein smoothing theory. This is the first
 example of such complex surfaces.

 The key ingredient of this paper is to find an elliptic surface $Y$
 equipped with a special \emph{bisection} ($=$ an irreducible curve
 on an elliptic surface whose intersection number with a fiber is $2$).
 Blowing-up $Y$ several times appropriately, we get a rational
 surface $Z$ which makes it possible to get such a complex surface.
 Once we have a right candidate $Z$ for $K^2=4$, the remaining
 argument is similar to that of $K^2=1,2,3$ case appeared
 in~\cite{LP07,PPS}. That is, by applying a rational blow-down surgery
 and $\mQ$-Gorenstein smoothing theory developed in~\cite{LP07} to
 $Z$, we obtain a minimal complex surface of general type with
 $p_g=0$ and $K^2 =4$. Then we show that the surface is simply
 connected. Since almost all the proofs are parallel to the case of
 the main construction in~\cite[\S3]{PPS}, we only explain how to
 construct such a minimal complex surface and we prove that the surface
 is simply connected. The main result of this paper is the following

\begin{theorem}
\label{thm-main}
 There exists a simply connected minimal complex surface of general type
 with $p_g=0$ and $K^2=4$.
\end{theorem}

\begin{remark}
 R\u{a}sdeaconu and S\c{u}vaina \cite{RS}
 proved that the complex surfaces constructed in~\cite{LP07,PPS} admit
 a K\"ahler-Einstein metric of negative scalar curvature.
 By applying their method to the complex surface constructed in this paper, one may prove that it also admits a K\"ahler-Einstein metric
 of negative scalar curvature; \S\ref{section:einstein}.
\end{remark}

\b

\section{Main construction}\label{section:main-construction}

\begin{figure}[hbtb]
 \centering \subfloat[Two generators]
 {\includegraphics[scale=0.7]{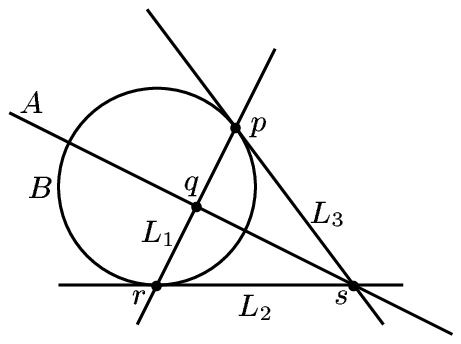}} 
 \centering \subfloat[A bisection]
 {\includegraphics[scale=0.7]{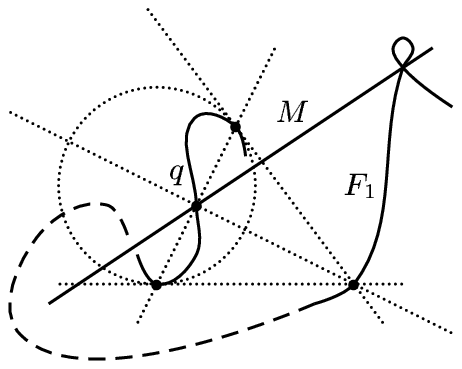}}
 \caption{A pencil of cubics}
 \label{figure:pencil}
\end{figure}

\begin{figure}[hbtb]
 \centering
 \includegraphics[scale=0.7]{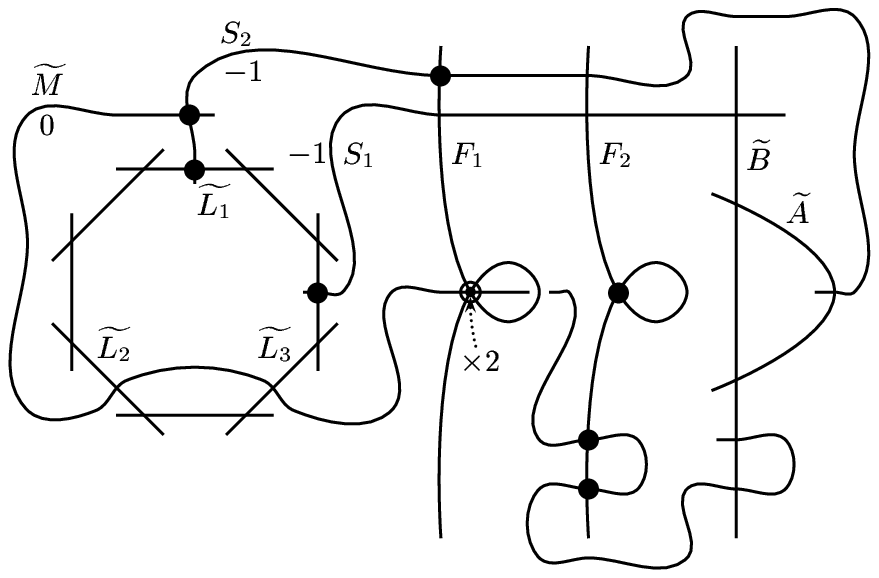}
 \caption{An elliptic fibration $Y$}
 \label{figure:Y}
 \vspace*{-1 em}
\end{figure}

 We start with a special elliptic fibration $Y:=\mP^2\sharp 9\overline{\mP}^2$
 which is used in the main construction of this paper. Let $L_1$,
 $L_2$, $L_3$ and $A$ be lines in $\mP^2$ and let $B$ be a smooth
 conic in $\mP^2$ intersecting as in Figure~\ref{figure:pencil}(a).
 We consider a pencil of cubics $\{\lambda(L_1+L_2+L_3) + \mu(A+B)
 \mid [\lambda:\mu] \in \mP^1 \}$ in $\mP^2$ generated by two cubic
 curves $L_1+L_2+L_3$ and $A+B$, which has $4$ base points, say, $p$,
 $q$, $r$ and $s$. In order to obtain an elliptic fibration over
 $\mP^1$ from the pencil, we blow up three times at $p$ and $r$,
 respectively, and twice at $s$, including infinitely near
 base-points at each point. We perform one further blowing-up at the
 base point $q$. By blowing-up totally nine times, we resolve all base points
 (including infinitely near base-points) of the pencil and we then
 get an elliptic fibration $Y=\mP^2\sharp 9\overline{\mP}^2$ over
 $\mP^1$ (Figure~\ref{figure:Y}).

 There are four sections of the elliptic fibration $Y$ corresponding to
 the four base points $p$, $q$, $r$, and $s$. Among these sections we
 use only two sections corresponding to $p$ and $q$, say $S_1$ and
 $S_2$ respectively, for the main construction. Furthermore, the
 elliptic fibration $Y$ has an $I_8$-singular fiber consisting of the
 proper transforms $\widetilde{L_i}$ of $L_i$ ($i=1,2,3$). Also $Y$
 has an $I_2$-singular fiber consisting of the proper transforms
 $\widetilde{A}$ and $\widetilde{B}$ of $A$ and $B$, respectively.
 According to the list of Persson~\cite{Pers}, we may assume that $Y$
 has only two more nodal singular fibers $F_1$ and $F_2$ by choosing
 generally $L_i$'s, $A$ and $B$ (Figure~\ref{figure:Y}). For example
 the pencil used in~\cite{PPS} works:
\begin{equation}
\label{equation:explicit-pencil}
 \{\lambda (y-\sqrt{3}x)(y+\sqrt{3}x)(2y-3z) + \mu x(x^2 + (y-2z)^2-z^2)
 \mid [\lambda:\mu] \in \mP^1 \}.
\end{equation}
 This pencil has singular fibers at
 $[\lambda:\mu]=[1:0], [0:1], [2:3\sqrt{3}]$ and $[2:-3\sqrt{3}]$.
 Furthermore, setting $F_1 = \{2(y-\sqrt{3}x)(y+\sqrt{3}x)(2y-3z) + 3\sqrt{3}
 x(x^2 + (y-2z)^2-z^2)= 0 \}$ and $F_2 = \{ 2(y-\sqrt{3}x)(y+\sqrt{3}x)(2y-3z)
 - 3\sqrt{3}x(x^2 + (y-2z)^2-z^2) = 0 \}$, one can easily check that $F_1$ and
$F_2$ are nodal cubic curves with one node at $[\sqrt{3}:0:-1]$ and
$[\sqrt{3}:0:1]$, respectively.

 Let $M$ be the line in $\mP^2$ passing through the point $q$ and the node
 of the nodal cubic curve $F_1$. The node of $F_1$ does not lie on
 any $L_i$'s, $A$, and $B$. Hence it satisfies that $M \neq L_1$, $M \neq A$,
 and $\widetilde{M}\cdot \widetilde{M} = 0$, where $\widetilde{M}$
 is the proper transform of $M$ in $Y$ (Figure~\ref{figure:pencil}(b)).
 We may assume further
 that $M$ does not pass through the node of the other nodal cubic
 curve $F_2$ by choosing generally $L_i$'s, $A$, and $B$. For
 example, the pencil in \eqref{equation:explicit-pencil} works: We
 have $q=[0:3:2]$; hence the line $M$ passing through $q$ and the
 node of $F_1$ is $\{s[0:3:2]+t[\sqrt{3}:0:-1] \mid [s:t] \in
 \mP^1\}$. It is obvious that the node $[\sqrt{3}:0:1]$ does not lie
 on the line $M$. Since $M$ meets every member in the pencil at three
 points, $\widetilde{M}$ is a bisection of the elliptic fibration $Y
 \to \mP^1$. Furthermore, since $q \in M$, the section $S_2$ meets
 $\widetilde{M}$ at one point (Figure~\ref{figure:Y}).

\begin{figure}[hbtb]
 \centering
 \includegraphics[scale=0.7]{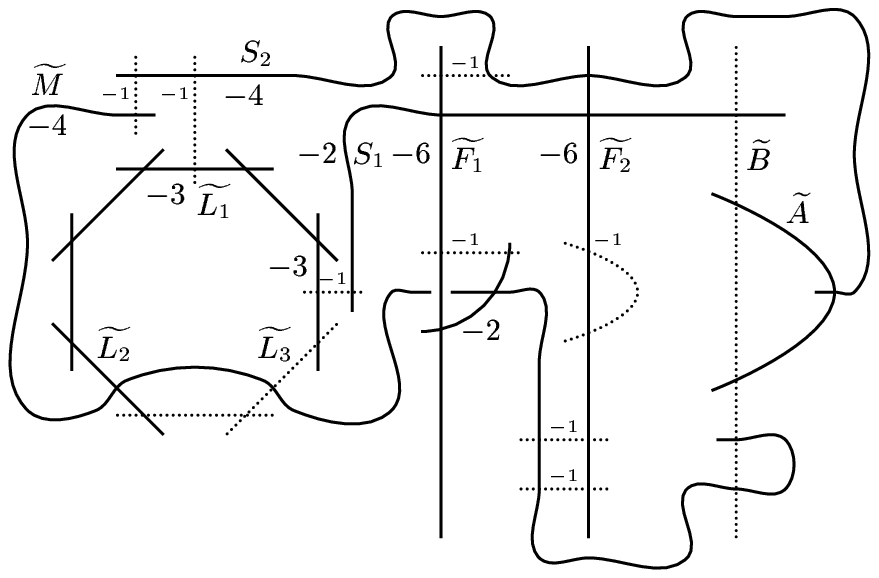}
 \caption{A rational surface $Z=Y\sharp 9{\overline \mP}^2$}
 \label{figure:Z}
 \vspace*{-1 em}
\end{figure}

 Next, by blowing-up nine times on $Y$, we construct a rational surface $Z$
 which contains a special configuration of linear chains of
 $\mP^1$'s. At first we blow up twice at the marked point
 $\bigodot$ on $F_1$. We then blow up seven times totally at the six
 marked points $\bullet$ on each fibers and at the intersection point
 $\bullet$ of $\widetilde{M}$ and $S_2$. We then get a rational
 surface $Z = Y \sharp 9\overline{\mP}^2$. We also denote by
 $\widetilde{F_i}$ ($i=1,2$) the proper transforms of $F_i$. Then
 there exists a linear chain of ${\mP}^1$'s in $Z$:
 \[C_{252,145} =\ntc{-2}{u_{13}}-\ntc{-4}{u_{12}}-\ntc{-6}{u_{11}}
 -\ntc{-2}{u_{10}}-\ntc{-6}{u_9}-\ntc{-2}{u_8}-\ntc{-4}{u_7}-\ntc{-2}{u_6}
 -\ntc{-2}{u_5}-\ntc{-2}{u_4}-\ntc{-3}{u_3}-\ntc{-2}{u_2}-\ntc{-3}{u_1}, \]
 which contains $\widetilde{A}$, $S_2$, $\widetilde{F_2}$, $S_1$,
 $\widetilde{F_1}$, $\widetilde{M}$, $\widetilde{L_2}$,
 $\widetilde{L_1}$, and $\widetilde{L_3}$, where $u_i$ represents an
 embedded rational curve (Figure~\ref{figure:Z}).

 Finally, by applying $\mQ$-Gorenstein smoothing theory to $Z$ as
 in~\cite{LP07, PPS}, we construct a minimal complex surface with
 $p_g=0$ and $K^2=4$. That is, we first contract the chain
 $C_{252,145}$ of $\mP^1$'s from
 $Z$ so that it produces a normal projective surface $X$ with one
 permissible singular point. And then, it follows by a similar technique
 in~\cite{LP07, PPS} that $X$ has a $\mQ$-Gorenstein smoothing. Let
 $X_t$ be a general fiber of the $\mQ$-Gorenstein smoothing of $X$.
 Since $X$ is a (singular) surface with $p_g=0$ and $K^2=4$, by
 applying general results of complex surface theory and
 $\mQ$-Gorenstein smoothing theory, one may conclude that a general
 fiber $X_t$ is a complex surface of general type with $p_g=0$ and
 $K^2=4$.

 The minimality of $X_t$ follows from the nefness of the canonical
 divisor $K_X$ of $X$. Let $f : Z \to X$ be the contraction of the
 chain $C_{252,145}$ of $\mP^1$'s from $Z$ to the singular surface
 $X$. By using a similar technique in~\cite{LP07, PPS}, it follows that
 the pullback $f^{\ast}{K_X}$ of the canonical divisor $K_X$ of $X$
 is effective and nef, hence $K_X$ is also nef, which shows the
 minimality of $X_t$.

 It remains to prove that $X_t$ is simply connected.

\begin{proposition}\label{proposition:simply-connected-I}
 $X_t$ is simply connected.
\end{proposition}

\begin{proof}
 Let $Z_{252}$ be a rational blow-down $4$-manifold obtained from $Z$
 by replacing the configuration $C_{252,145}$ with the corresponding
 rational ball $B_{252,145}$. Since a general fiber $X_t$ of a
 $\mQ$-Gorenstein smoothing of $X$ is diffeomorphic to the rational
 blow-down $4$-manifold $Z_{252}$, it suffices to show that
 $Z_{252}$ is simply connected. We decompose the surface $Z$ into
 $Z= Z_0 \cup C_{252, 145}$. Then we have $Z_{252} = Z_0 \cup B_{252, 145}$.
 Furthermore,
 since $\pi_{1}(\partial B_{252,145}) \rightarrow \pi_{1}(B_{252,145})$
 is surjective, by Van-Kampen theorem,
 it suffices to show that $\pi_{1}(Z_0) =1$.

 Let $\alpha_i$ be a normal circle of $u_i$. First, note that $Z$ and
 the configuration $C_{252, 145}$ are all simply connected. Hence,
 applying Van-Kampen theorem on $Z$, we get
\begin{equation}
\label{equation:VanKampen}
 1 = \pi_{1}(Z_0)/\langle N_{i_{*}(\alpha_1)} \rangle,
\end{equation}
 where $i_{*}$ is the induced homomorphism by the inclusion
 $i : \partial C_{252, 145} \rightarrow Z_0$.

\begin{figure}[hbtb]
 \centering
 \includegraphics[scale=0.6]{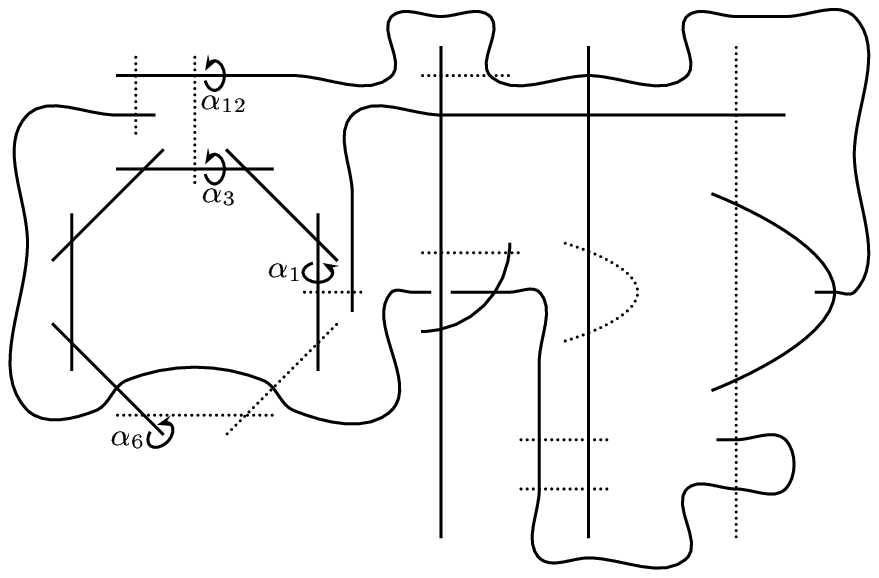}
 \caption{Normal circles}
 \label{figure:circle}
 \vspace*{-1 em}
\end{figure}

 We denote $a \sim b$ if $a$ and $b$ are conjugate to each other
 in $\pi_{1}(Z_0)$. From Figure~\ref{figure:circle},
 one can easily show that $1=i_{*}(\alpha_6) \sim i_{*}(\alpha_1)^{26}$,
 i.e. $i_{*}(\alpha_1)^{26}=1$, and
 $i_{*}(\alpha_1)^5 \sim i_{*}(\alpha_3) \sim i_{*}(\alpha_{12}) \sim
 i_{*}(\alpha_1)^{9574}.$
 Since $9574 \equiv 6\pmod{26}$, we have $i_{*}(\alpha_1)^5 \sim
 i_{*}(\alpha_1)^6$.
 Hence $i_{*}(\alpha_1)^{5\cdot13} \sim i_{*}(\alpha_1)^{26} = 1$,
 which implies that $i_{*}(\alpha_1)^{5 \cdot 13} = 1$. Since
 $\alpha_1^{5 \cdot 13}$ is also a generator of $\pi_1(\partial
 C_{252, 145})$,
 we have $i_{*}(\alpha_1)=1$.
 Therefore $\pi_{1}(Z_0) =1$ by \eqref{equation:VanKampen}.
\end{proof}

\b

\section{More Examples}\label{section:more}

 In this section we describe another rational surface $Z$ which makes it
 possible to get a simply connected surface of general type with $p_g=0$
 and $K^2=4$.

\subsection*{Construction}

\begin{figure}[hbtb]
 \centering \subfloat[Two generators]
 {\includegraphics[scale=0.7]{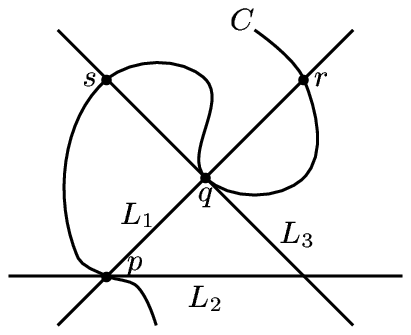}} 
 \centering \subfloat[Two lines $T$ and $M$]
 {\includegraphics[scale=0.7]{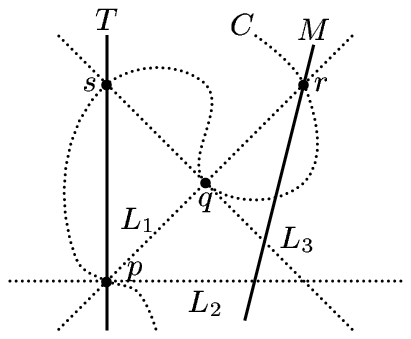}}
 \caption{A pencil of cubics}
 \label{figure:pencil2}
\end{figure}

 Let $C$ be a smooth cubic curve in $\mP^2$ and $p$ its inflection
 point. Let $L_1$ be a line passing through $p$ which intersects $C$ at two
 more different points $q$ and $r$. Let $L_2$ be the tangent line to
 $C$ at $p$ and $L_3$ the tangent line to $C$ at one of intersection
 points of $L_1$ and $C$, say $q$. Let $s$ be the other intersection
 point of $L_3$ and $C$ (Figure~\ref{figure:pencil2}(a)).
 We consider a pencil of cubics $\{\lambda(L_1+L_2+L_3) + \mu C \mid
 [\lambda:\mu] \in \mP^1 \}$ in $\mP^2$ generated by two cubic curves
 $L_1+L_2+L_3$ and $C$. According to \cite{Pers}, if we choose a general $C$,
 we may assume that the pencil of cubics contains four nodal singular
 curves. Let $T$ be a line joining $p$ and $s$ and $M$ a line through
 $r$ and the node of a nodal singular member of the pencil of cubics.
 We may assume that $M$ does not pass through the other nodes
(Figure~\ref{figure:pencil2}(b)).

\begin{figure}[hbtb]
 \centering
 \includegraphics[scale=0.7]{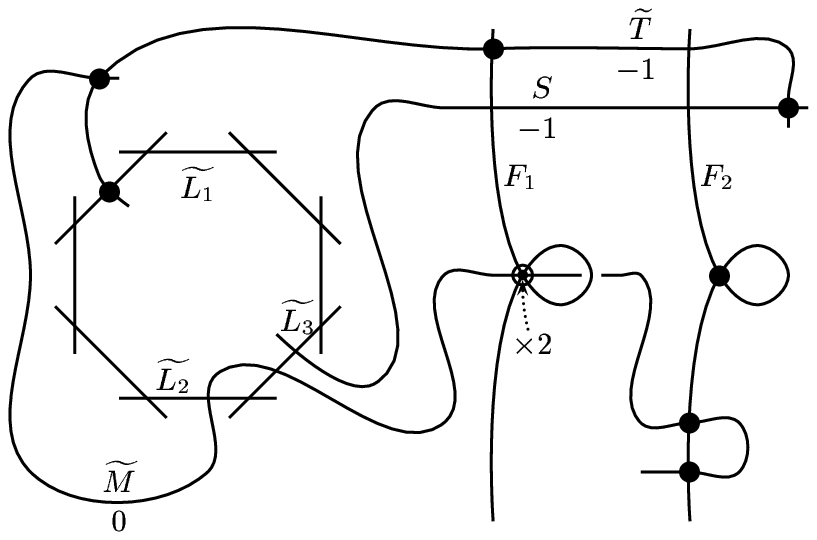}
 \caption{An elliptic fibration $Y$}
 \label{figure:Y2}
 \vspace*{-1 em}
\end{figure}

 In order to obtain an elliptic fibration over $\mP^1$ from the pencil above,
 we blow up totally $9$ times at the base points of the pencil of
 cubics including infinitely near base-points at each base point. We
 then get an elliptic fibration $Y=\mP^2\sharp 9\overline{\mP}^2$
 over $\mP^1$ (Figure~\ref{figure:Y2}). Note that the proper
 transform $\widetilde{T}$ of $T$ is a section of $Y$ and the proper
 transform $\widetilde{M}$ of $M$ is a bisection of $Y$
(Figure~\ref{figure:Y2}). Here the section $S$ in $Y$ is an
 exceptional curve induced by the blowing-up at the point $s$.

\begin{figure}[hbtb]
 \centering
 \includegraphics[scale=0.7]{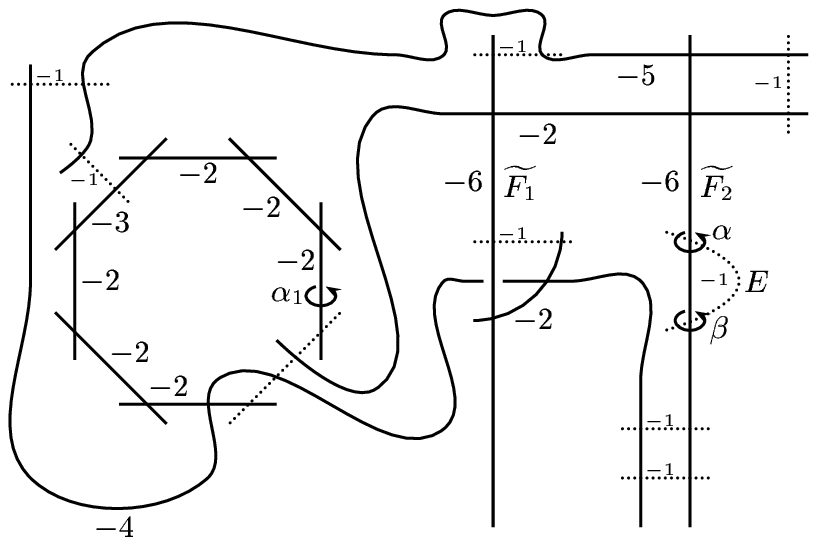}
 \caption{A rational surface $Z=Y\sharp 9{\overline \mP}^2$}
 \label{figure:Z2}
 \vspace*{-1 em}
\end{figure}

 We blow up $7$ times at the marked points $\bullet$ on $Y$ and blow up
 two more times at the marked point $\bigodot$ on $Y$. We finally
 obtain a rational surface $Z=Y\sharp 9{\overline \mP}^2$ which contains
 the following linear chain of $\mP^1$'s (Figure~\ref{figure:Z2}):
\[C_{183,38} =\ntc{-5}{u_{13}}-\ntc{-6}{u_{12}}-\ntc{-2}{u_{11}}
 -\ntc{-6}{u_{10}}-\ntc{-2}{u_9}-\ntc{-4}{u_8}-\ntc{-2}{u_7}-\ntc{-2}{u_6}
 -\ntc{-2}{u_5}-\ntc{-3}{u_4}-\ntc{-2}{u_3}-\ntc{-2}{u_2}-\ntc{-2}{u_1}.\]

 Finally, by applying $\mQ$-Gorenstein smoothing theory to $Z$
 as in~\cite{LP07, PPS}, we are able to construct a minimal complex
 surface with $p_g=0$ and $K^2=4$, say $X_t$, which is a general fiber
 of a $\mQ$-Gorenstein smoothing of $X$.

\begin{proposition}
The complex surface $X_t$ is simply connected.
\end{proposition}

\begin{proof}
 Let us decompose the surface $Z=Y\sharp 9{\overline \mP}^2$ into
 $Z= Z_0 \cup C_{183,38}$. Then, as in the proof of
 Proposition~\ref{proposition:simply-connected-I},
 it is enough to show that $\pi_{1}(Z_0) =1$.

 Let $E$ be an exceptional curve intersecting $\widetilde{F_2}$ at two points.
 The intersection of a boundary of a tubular neighborhood of $\widetilde{F_2}$
 and $E$ consists of two normal circles of $\widetilde{F_2}$, say $\alpha$
 and $\beta$, which are contained in $Z_0$.
 We choose a point $x_0 \in \alpha$ as a base point for the homotopy group
 of $Z_0$. Let $x_1 \in \beta$ be any point.

 Since $\widetilde{F_2}$ and $E$ intersect positively at each intersection
 point, $\alpha$ and $\beta$ have the same orientation induced by
 the orientation of the exceptional curve $E$. Therefore, as circles on the
 punctured sphere $E \setminus C_{183,38}$, they are the boundaries of
 the cylinder $E \setminus C_{183,38}$ and,
 furthermore, they have the opposite orientation in the
 cylinder $E \setminus C_{183,38}$. Let $i_{\ast}$ be the induced
 homomorphism by the inclusion $i : \partial C_{183,38} \to Z_0$.
 Then we have
\begin{equation}\label{equation:a=b^-1}
 [i_{\ast}(\alpha)] = [\lambda \cdot i_{\ast}(\beta)^{-1} \cdot \lambda^{-1}]
 \quad \text{in $\pi_1(Z_0, x_0)$},
\end{equation}
 where $\lambda$ is a path connecting $x_0$ and $x_1$ which lies on $E$.

 On the other hand, since $\alpha$ and $\beta$ are normal circles of
 $\widetilde{F_2}$, we also have
\begin{equation}\label{equation:a=b}
 [i_{\ast}(\alpha)] = [\mu \cdot i_{\ast}(\beta) \cdot \mu^{-1}]
 \quad \text{in $\pi_1(Z_0, x_0)$}
\end{equation}
 where $\mu$ is a path connecting $x_0$ and $x_1$ which is contained in
 the boundary of a tubular neighborhood of $\widetilde{F_2}$.
 Note that we may choose $\lambda$ and $\mu$ so that they are
 homotopically equivalent. Therefore it follows by
\eqref{equation:a=b^-1} and \eqref{equation:a=b} that
\begin{equation}\label{equation:a=1}
 [i_{\ast}(\alpha)^2] = 1 \quad \text{in $\pi_1(Z_0, x_0)$}.
\end{equation}

 It is not difficult to show that $i_{\ast}(\alpha)^2$ is conjugate to
 $i_{\ast}(\alpha_1)^{2552}$, where $\alpha_1$ is a generator of
 $\pi_1(\partial Z_0 =L(183^2, -6953), x_0)=\mZ_{183^2}$.
 Since $2552 = 8 \cdot 11 \cdot 29$ is
 relatively prime to $183^2 = (3 \cdot 61)^2$, it implies that
 $\alpha^2$ is also a generator of $\pi_1(\partial Z_0)$.
 By applying Van-Kampen theorem on $Z$, we get
\begin{equation*}
 1 = \pi_1(Z_0, x_0) / \langle N_{i_{\ast}(\alpha)^2} \rangle.
\end{equation*}
 Therefore $\pi_1(Z_0, x_0) = 1$ by \eqref{equation:a=1}.
\end{proof}

\begin{remark}
 1. One can find more examples of simply connected surface of general type
 with $p_g=0$ and $K^2=4$ using different configurations.
 For example, using an elliptic fibration on $E(1)$ with one
 $I_7$-singular fiber, one $I_2$-singular fiber and two nodal
 fibers, we can find the following linear chain of $\mP^1$'s in
 $E(1)\sharp 9{\overline \mP}^2$:
 \[C_{252,145} =\nc{-2}-\nc{-4}-\nc{-6}
 -\nc{-2}-\nc{-6}-\nc{-2}-\nc{-4}-\nc{-2}
 -\nc{-2}-\nc{-2}-\nc{-3}-\nc{-2}-\nc{-3}.\]
 It is a very intriguing question whether all these configurations above
 produce the same deformation equivalent type of simply connected surfaces
 with $p_g=0$ and $K^2=4$. We leave this problem for future research.\\
 2. It is also a natural question whether one can find an appropriate
 configuration in a rational surface which produces a surface of general
 type with $p_g=0$ and $K^2 \geq 5$.
 Note that the basic scheme used in this paper as well as
 in~\cite{LP07, PPS} is the following: We chose a delicate configuration
 in a certain rational surface $Z$ so that its induced singular surface $X$
 obtained by contracting linear chains of curves in $Z$ satisfies the
 cohomology condition $H^2(T^0_X)=0$, which guarantees automatically
 the existence of a $\mQ$-Gorenstein smoothing of $X$.
 In this respect, it seems impossible to find a configuration satisfying
 $H^2(T^0_X)=0$ for $K^2 \geq 5$.
 But, without hypothesis $H^2(T^0_X) =0$,
 it might still be a chance to find a configuration for $K^2 \geq 5$.
 Of course, if such a configuration exists, it will be another problem
 to determine whether the induced singular surface $X$ admits
 a $\mQ$-Gorenstein smoothing or not.
\end{remark}

\b

\section{Einstein metrics on $\mC\mP^2 \sharp 5\overline{\mC\mP}^2$}
\label{section:einstein}

 In this section we show that the complex surface $X_t$ constructed
 in the main construction admits a K\"ahler-Einstein metric of negative scalar
 curvature, which implies the following theorem.

\begin{theorem}
\label{thm-einstein}
 The topological $4$-manifold $\mC\mP^2\sharp5\overline{\mC\mP}^2$
 has a smooth structure which admits an Einstein metric with negative
 scalar curvature.
\end{theorem}

 Recently R\u{a}sdeaconu and S\c{u}vaina \cite{RS} proved the
 existence of a smooth structure on each of the topological
 $4$-manifolds $\mC\mP^2\sharp k\overline{\mC\mP}^2$, for $k = 6, 7$,
 which has an Einstein metric of negative scalar curvature.
 By applying their method on a surface $X_t$ constructed in Section 2,
 we can easily prove the existence of a
 K\"ahler-Einstein metric on $X_t$ with negative scalar curvature.
 We explain it in a detail in the rest of this section.

 First, note that there is a criterion for the existence of a K\"ahler-Einstein
 metric on a compact complex $4$-manifold with $c_1(M) < 0$, which is found
 independently by Aubin \cite{Aubin} and Yau \cite{Yau}:

\begin{theorem}[Aubin-Yau]
 A compact complex $4$-manifold $(M, J)$ admits a compatible
 K\"ahler-Einstein metric with negative scalar curvature if and only
 if its canonical line bundle $K_M$ is ample. When such a metric exists,
 it is unique, up to an overall multiplicative constant.
\end{theorem}

\begin{proof}[Proof of Theorem~\ref{thm-einstein}]
 Based on the idea \cite{RS}, we show that the surface $X_t$ has an
 ample canonical bundle. Then it follows from the criterion of
 Aubin-Yau that there exists a K\"ahler-Einstein metric on $X_t$ of
 negative scalar curvature.

 As we showed in the main construction, the pullback $f^{\ast}{K_X}$
 of the canonical divisor $X$ onto the rational surface $Z$ is
 effective and nef; hence $K_X$ is also nef. Let $E_1, \dotsc, E_8$
 be the ($-1$)-curves on the rational surface $Z$ and set
\[ C_{252,145} =\ntc{-2}{G_{13}}-\ntc{-4}{G_{12}}-\ntc{-6}{G_{11}}-
   \ntc{-2}{G_{10}}-\ntc{-6}{G_9}-\ntc{-2}{G_8}-\ntc{-4}{G_7}-\ntc{-2}{G_6}-
   \ntc{-2}{G_5}-\ntc{-2}{G_4}-\ntc{-3}{G_3}-\ntc{-2}{G_2}-\ntc{-3}{G_1}, \]
 Then one may write
\[ f^{\ast}{K_X} \equiv_{\mQ} \sum_{i=1}^{8}{a_i E_i} +
   \sum_{j=1}^{13}{b_j G_j} \]
 for some rational numbers $a_i, b_j \ge 0$.

 We first show that $K_X$ is ample. Suppose on the contrary that
 $K_X$ is not ample. Since $K_X$ is already nef and $K_X^2 = 4 > 0$,
 according to the Nakai-Moishezon criterion, there exists an
 irreducible curve $C \subset X$ such that $(K_X \cdot C) = 0$.
 Let $\overline{C} \subset Z$ be  the proper transform of $C$.
 Then we have
\[(K_X \cdot C) = (f^{\ast}{K_X} \cdot f^{\ast}{C}) =
  (f^{\ast}{K_X} \cdot \overline{C}) =
  \sum_{i=1}^{8}{a_i (E_i \cdot \overline{C})} +
  \sum_{j=1}^{13}{b_j (G_j \cdot \overline{C})} = 0. \]

 Since $G_j$'s are irreducible components of the exceptional divisors
 of $f$, it is obvious that $(G_j \cdot \overline{C}) \ge 0$
 ($j=1, \dotsc, 13$) with equality if and only if $C$ does not pass
 through the singular point of $X$. Hence it follows that
\[ \sum_{i=1}^{8}{a_i (E_i \cdot \overline{C})} \le 0. \]
 Then either $(E_{i_0} \cdot \overline{C}) < 0$ for some $i_0$, or
 $(E_i \cdot \overline{C})= 0$ for all $i = 1, \dotsc, 8$ and
 $(G_j \cdot \overline{C}) = 0$ for all $j=1, \dotsc, 13$. In the first
 case $\overline{C}$ must coincide with $E_{i_0}$. However, by using
 a similar technique in~\cite{LP07, PPS}, one may show that
 $(f^{\ast}{K_X} \cdot E_i) > 0$ for all $i = 1, \dotsc, 8$, which is a
 contradiction to our assumption $(K_X \cdot \overline{C}) = 0$.
 Therefore we have $(E_i \cdot \overline{C})= 0$ for all $i = 1, \dotsc, 8$
 and $(G_j \cdot \overline{C}) = 0$ for all $j=1, \dotsc, 13$.
 On the other hand, note that the Poincar\'e duals of the irreducible components
 $G_j$'s and of the ($-1$)-curves $E_i$'s generate $H_2(Z,\mQ)$;
 hence $\overline{C}$ must be numerically trivial on $Z$. Then, for
 any ample divisor $H$ on $X$, we have
\[ 0 = (\overline{C} \cdot f^{\ast}{H}) = (f^{\ast}{C} \cdot f^{\ast}{H})
     = (C \cdot H), \]
 which is again a contradiction. Therefore $K_X$ is ample.

 Note that ampleness is an open property \cite{KM}. So the canonical
 divisor $K_{X_t}$ of a general fiber $X_t$ of $\mQ$-Gorenstein
 smoothing is automatically ample. Therefore, by Aubin-Yau's
 criterion, $X_t$ has a K\"ahler-Einstein metric of negative scalar
 curvature.
\end{proof}

\m

\subsubsection*{Acknowledgements}
 The authors would like to thank Yongnam Lee for helpful
 communications during the course of this work.
 Jongil Park was supported by the Korea Research Foundation Grant funded
 by the Korean Government (KRF-2007-314-C00024) and he also holds a
 joint appointment in the Research Institute of Mathematics,
 Seoul National University.
 Dongsoo Shin was supported by Korea Research Foundation Grant funded
 by the Korean Government (KRF-2005-070-C00005).

\b
\b

\providecommand{\bysame}{\leavevmode\hbox
 to3em{\hrulefill}\thinspace}


\begin{thebibliography}{999}
 \bibitem[1]{Aubin} T. Aubin, {\it {\'E}quations du type Monge-Amp{\`e}re sur
          les vari\'et\'es k\"ahleriennes compactes}, C. R. Acad. Sci. Paris
          Sr. A-B \textbf{283} (1976), no. 3, A119--A121.

 \bibitem[2]{B} R. Barlow, {\it A simply connected surface of general type
                with $p_g =0$}, Invent. Math. \textbf{79} (1984), 293--301.

 \bibitem[3]{BHPV} W. Barth, K. Hulek, C. Peters and A. Van de Ven,
                  {\it Compact complex surfaces}, 2nd ed. Springer-Verlag, 2004.

 \bibitem[4]{KM} J. Koll\'ar, S. Mori, {\it Birational geometry of algebraic
                  varieties}, Cambridge Tracts in Mathematics, \textbf{134},
                  Cambridge University Press, Cambridge, 1998.

 \bibitem[5]{LP07} Y. Lee and J. Park, {\it A simply connected surface of
                 general type with $p_g=0$ and $K^2=2$},
                 Invent. Math. \textbf{170} (2007), 483--505.

 \bibitem[6]{P05} J. Park, {\it Simply connected symplectic 4-manifolds
                  with $b_2^+=1$ and $c_1^2=2$},
                  Invent. Math. \textbf{159} (2005), 657--667.

 \bibitem[7]{PPS} H. Park, J. Park and D. Shin, {\it A simply connected
                  surface of general type with $p_g=0$ and $K^2=3$},
                  arXiv:0708.0273.

 \bibitem[8]{Pers} U. Persson, {\it Configuration of Kodaira fibers on
                  rational elliptic surfaces},
                  Math. Z. \textbf{205} (1990), 1--47.

 \bibitem[9]{RS} R. R\u{a}sdeaconu and I. S\c{u}vaina, {\it Smooth structures
           and einstein metrics on $\mC\mP^2\sharp5, 6, 7\overline{\mC\mP}^2$},
           arXiv:0806.1424.

 \bibitem[10]{Yau} S.-T. Yau, {\it Calabi¡¯s conjecture and some new results
            in algebraic geometry}, Proc. Nat. Acad. USA \textbf{74} (1997),
            1789--1799.

\end{thebibliography}
\end{document}